\documentclass[12pt]{article}
\usepackage{amsmath}
\usepackage{amssymb}
\usepackage{amsthm}
\usepackage{graphicx}
\usepackage[usenames]{color}

\usepackage{hyperref}

\newtheorem{thm}{Theorem}[section]

\newtheorem{rem}[thm]{Remark}

\newtheorem{definition}[thm]{Definition}

\renewcommand{\qed}{\ifmmode$\Box$\else{\unskip\nobreak\hfil
	\penalty50\hskip1em\null\nobreak\hfil$\Box$
	\parfillskip=0pt\finalhyphendemerits=0\endgraf}\fi}

\newcommand{\eps}{\varepsilon}
\newcommand{\bfu}{{\boldsymbol u}}
\newcommand{\bfv}{{\boldsymbol v}}
\newcommand{\bfz}{{\boldsymbol z}}

\newcommand{\bfg}{{\boldsymbol g}}
\newcommand{\bfG}{{\boldsymbol G}}
\newcommand{\bfW}{{\boldsymbol W}}
\newcommand{\bfH}{{\boldsymbol H}}
\newcommand{\bfV}{{\boldsymbol V}}

\makeatletter
\def\blfootnote{\xdef\@thefnmark{}\@footnotetext}
\makeatother

\begin{document}

\begin{center}
{\bf \large Variational Inequalities of Navier--Stokes Type\\[0.25cm]
            with Time Dependent Constraints}
\vspace{0.5cm}

\blfootnote{\copyright 2016. This manuscript version is made available under the \href{http://creativecommons.org/licenses/by-nc-nd/4.0/ }{CC-BY-NC-ND 4.0}  license} 

Maria Gokieli, Nobuyuki Kenmochi and
Marek Niezg\'odka

Interdisciplinary Centre of Mathematical and Computational Modelling,

       University of Warsaw,

       Pawi\'nskiego 5a, 02-106 Warsaw, Poland
\end{center}
\vspace{1.5cm}

\noindent
{\bf Abstract.} We consider a class of parabolic variational
inequalities with time dependent obstacle of the form 
$|{\boldsymbol u}(x,t)| \le p(x,t)$, where ${\boldsymbol u}$ is the velocity
field of a fluid governed by the Navier--Stokes variational inequality. 
The obstacle function $p=p(x,t)$ imposed on ${\boldsymbol u}$
consists of three parts which are respectively the degenerate part $p(x,t)=0$, 
the finitely positive part $0< p(x,t) <\infty$ and singular part
$p(x,t)=\infty$. In this paper, we shall propose a sequence of approximate 
obstacle problems with everywhere finitely positive obstacles and prove an
existence result for the original problem by discussing the convergence of
the approximate problems. The crucial step is to handle the nonlinear
convection term. In this paper we propose a new approach to it.

\renewcommand{\large}{\normalsize}
\renewcommand{\Large}{\normalsize}
\renewcommand{\LARGE}{\normalsize}

\section{Introduction}\label{sec:intro}

In real problems, we find many dynamical processes which occur in fluids or in consequence of a fluid flow. 
Their mathematical models include then a hydrodynamic equation, 
typically Stokes or Navier--Stokes, coupled with some other
evolution systems, such as heat 
convection [14, 16], phase transitions [1] or biofilm growth [11, 23].
These couplings may have the form of~transport or advection, 
but they may also mean some evolution of the domain in which 
the flow takes place. To give just one example
of phenomenon of importance to medicine and in which 
both types of couplings appear at the same time scale, 
let us have a look on the mentioned biomass growth. 
In a fluid transporting some living organisms and some appropriate nutrient,
some of these organisms can stick to the boundary of the fluid flow's domain (e.g.~blood vessels walls) and then aggregate, in which way they gradually restrict the domain available for the flow, forming a geometrical obstacle to~it.

Mathematical analysis of  models for such systems  seems not easy, 
since the theory on
partial differential systems coupled with equations, or variational inequalities, of the Navier--Stokes type has not been completely established.

In this paper, we address the problem of a Navier--Stokes flow
constrained by some evolving in time obstacle. 
We model the obstacle as a non-negative 
function $p$, depending on the space and time variable, 
which is a bound imposed {\it a priori} on the velocity of the flow. 
The Navier--Stokes equation becomes then naturally 
a variational inequality.
We allow the constraint to disappear ($p=\infty$, free flow), 
to be a total obstacle ($p=0$, no flow)
or only partial ($0<p<\infty$). 
We assume that $p$ is continuous. 
Our main result is Theorem~1.1 below, stating existence
and some regularity of solution to this problem.

This kind of parabolic obstacle problem would 
be useful for mathematical modelling of
various nonlinear problems in hydrodynamic 
fluids, see e.g.~[2, 9, 10, 11, 18, 19, 21]. 
As far as variational inequalities of Navier--Stokes 
type are concerned, see 
e.g.~[3--6, 24, 25] for a constant in time constraint,  
and [13], where the constraint can be time and space dependent.  
However, even this last case did not
allow the "free flow" and "no-flow" regions, 
i.e.~the obstacle function $p$ had to be finite 
and bounded from below by a positive constant ---
a serious limitation of the model that
we overcome in the present work. It is clear that especially 
allowing the "total obstacle" case, 
i.e.~having regions where $p=0$, 
is essential from the point of view of modelling; 
it is also the main challenge for the mathematical 
analysis that we are presenting.

For basic studies on Navier--Stokes equations and phase transitions, 
we refer to [27] and [8], respectively. 
Our formulation of the Navier--Stokes inequality 
arising from the obstacle, that we state in~Definition 1.1 below,
is analogous to these appearing in [3--6, 24, 25].
For its analysis, we will use the theory of subdifferentials  
contained in~[7, 17, 22, 28]. This will be exposed in Section~2.
\\

Let us set the basic functional framework and explicit the assumptions so as to formulate the main result.  
Let 
$\Omega$ be a bounded domain in ${\bf R}^3$ with smooth boundary
$\Gamma:=\partial \Omega$, $Q:=\Omega \times (0,T)$, $0<T<\infty$ and
$\Sigma:=\Gamma \times (0,T)$, and denote by $|\cdot|_X$ the norm in various
function spaces $X$ built on $\Omega$ as well as by $\|\cdot \|_Y$ for
function spaces $Y$ on $\Omega\times (0,T)$. Also, consider the usual solenoidal function spaces:
\begin{eqnarray*}
 && {\boldsymbol {\mathcal D}}_\sigma(\Omega):=
      \{{\boldsymbol v}=(v^{(1)},v^{(2)},v^{(3)})
       \in {\cal D}(\Omega)^3~|~
  {\rm div}\ {\boldsymbol v}=0~{\rm in~}\Omega\},\\
 &&{\boldsymbol H}_\sigma(\Omega):= {\rm the~ closure~ of~}{\boldsymbol {\mathcal D}}_\sigma(\Omega)~
    {\rm in}~ L^2(\Omega)^3,~{\rm with~norm~}|\cdot|_{0,2},\\
 &&{\boldsymbol V}_\sigma(\Omega):= {\rm the~ closure~ of~}{\boldsymbol {\mathcal D}}_\sigma(\Omega)~
    {\rm in~} H^1_0(\Omega)^3,~{\rm with~norm~} |\cdot|_{1,2},\\
 &&{\boldsymbol W}_\sigma(\Omega):={\rm the~ closure~ of~}
   {\boldsymbol {\mathcal D}}_\sigma(\Omega) ~{\rm in~}W^{1,4}_0(\Omega)^3,
      ~{\rm with~norm~}|\cdot|_{1,4};
\end{eqnarray*}
in these spaces the norms are given as usual by
 $$ |{\boldsymbol v}|_{0,2}:=\left\{ \sum_{k=1}^3 \int_\Omega |v^{(k)}|^2dx \right\}
   ^{\frac 12},~~ |{\boldsymbol v}|_{1,2}:=\left\{ \sum_{k=1}^3 \int_\Omega 
     |\nabla v^{(k)}|^2dx \right\}^{\frac 12}$$
and 
 $$ |{\boldsymbol v}|_{1,4}:=\left\{ \sum_{k=1}^3 \int_\Omega |\nabla v^{(k)}|^4dx 
   \right\}^{\frac 14}.$$   
For simplicity we denote 
the dual spaces of ${\boldsymbol V}_\sigma(\Omega)$ and 
${\boldsymbol W}_\sigma(\Omega)$ 
by ${\boldsymbol V}^*_\sigma(\Omega)$ and ${\boldsymbol W}^*_\sigma(\Omega)$, 
respectively, which are equipped with their dual norms $|\cdot|_{-1,2}$ and
$|\cdot|_{-1,\frac 43}$. 
Also, we denote the inner product in 
${\boldsymbol H}_\sigma(\Omega)$ by $(\cdot,\cdot)_\sigma$ and the duality 
between 
${\boldsymbol V}^*_\sigma(\Omega)$ and ${\boldsymbol V}_\sigma(\Omega)$ by 
$\langle \cdot,\cdot \rangle_\sigma$, namely for 
${\boldsymbol v}_i=(v^{(1)}_i,v^{(2)}_i,v^{(3)}_i),~i=1,2,$ 
$$ ({\boldsymbol v}_1,{\boldsymbol v}_2)_\sigma :=\sum_{k=1}^3\int_\Omega 
   v_1^{(k)}v_2^{(k)}dx, \qquad
\langle {\boldsymbol v}_1,{\boldsymbol v}_2\rangle_\sigma
    =\sum_{k=1}^3\int_\Omega \nabla v_1^{(k)}\cdot \nabla v_2^{(k)}dx. 
$$
Then, by identifying the dual of ${\boldsymbol H}_\sigma(\Omega)$ with itself, 
we have:
 $$ {\boldsymbol V}_\sigma(\Omega)
 \hookrightarrow
 {\boldsymbol H}_\sigma(\Omega)
 \hookrightarrow {\boldsymbol V}^*_\sigma(\Omega), 
 \quad
 {\boldsymbol W}_\sigma(\Omega)\hookrightarrow
  C(\overline{\Omega})^3;$$
and all these embeddings are compact.\vspace{0.5cm}

We are given a non-negative function $p=p(x,t)$ on $\overline Q$ as an
obstacle function such that
$ 0\le p(x,t) \le \infty$ for all $(x,t) \in \overline Q$ and $p$ is 
continuous from $\overline Q$ into $[0,\infty]$, namely,
 $$ \left\{
    \begin{array}{l}
 {\rm the~set~}\overline Q(p=\infty):=\{(x,t)\in \overline Q~|~p(x,t)=\infty\}
  ~{\rm is~closed~in~}
    \overline Q,\\[0.18cm] 
     \forall \kappa \in (0,\infty),~p~{\rm is~continuous~on~} 
      \overline Q(p\le \kappa):= \{(x,t)\in \overline Q~|~p(x,t) \le \kappa\},
      \\[0.18cm]
     \forall M \in (0,\infty), {\rm there~is~an~open~set~}U_M
     ~{\rm containing~}
      \overline Q(p=\infty) \\
     {\rm ~~~such ~that~} p \geq M~{\rm on~}U_M\cap \overline Q.
    \end{array} 
    \right. \eqno{(1.1)} $$
It is easily seen that (1.1) is equivalent to the  continuity on $\overline Q$ in the usual sense, of the function
  $$\alpha(x,t):=\left \{
        \begin{array}{l}
       \displaystyle{  \frac {p(x,t)}{1+p(x,t)},
     ~~~{\rm if~}0\le p(x,t) < \infty,}\\[0.5cm]
         1.~~~~~~~~~~~~~~~{\rm if~}p(x,t)=\infty,
        \end{array}
     \right. $$

We are now ready to define the solution of our obstacle problem.

\vspace{0.5cm}
\noindent
{\bf Definition 1.1.} {\it 
For given data
 $$  \nu >0~(\text{constant}),
     ~ {\boldsymbol g}\in L^2(0,T; {\boldsymbol H}_\sigma(\Omega)),~
   {\boldsymbol u}_0 \in {\boldsymbol H}_\sigma(\Omega),
 $$
our obstacle problem $P(p;{\boldsymbol g},{\boldsymbol u_0})$ 
is to find a solution 
${\boldsymbol u}:=(u^{(1)},u^{(2)},u^{(3)})$ from $[0,T]$ into 
${\boldsymbol H}_\sigma(\Omega)$ satisfying the following (i) and (ii):
\begin{description}
\item{(i)} ${\boldsymbol u}(0)
={\boldsymbol u}_0$ in~${\boldsymbol H}_\sigma(\Omega)$, and
$t\mapsto ({\boldsymbol u}(t), {\boldsymbol v}(t))_\sigma$ is of  
bounded variation on
$[0,T]$ for all ${\boldsymbol v} \in \boldsymbol{\cal K}(p)$, 
where
 $$ \boldsymbol{\cal K}(p):=\left\{{\boldsymbol v} \in C^1([0,T];{\boldsymbol W}_\sigma
(\Omega))~\left|
     \begin{array}{l}
      |{\boldsymbol v}|\le p \text{ on~}Q,~
      {\rm supp}({\boldsymbol v}) \text{ is~compact}\\
      \text{in~}\{(x,t)\in \Omega\times [0,T]~|~p(x,t)>0\}
     \end{array}
      \right. \right \}, $$
 and ${\rm supp}({\bfv})$ denotes the support of $\bfv$, 
 \item{(ii)} ${\boldsymbol u}: [0,T]\to {\bfH}_\sigma(\Omega),~\sup_{t \in [0,T]}|\bfu(t)|
_{0,2}< \infty$, ${\bfu}\in L^2(0,T;{\bfV}_\sigma(\Omega))$
and 
   $$|{\boldsymbol u}(x,t)| \le p(x,t)~\text{ a.e.~}x \in~\Omega,~\forall
     t \in [0,T],$$
   $$\int_0^t({\boldsymbol v}'(\tau), {\boldsymbol u}(\tau)-{\boldsymbol v}(\tau))_\sigma d\tau
    +\nu \int_0^t \langle {\boldsymbol u}(\tau), {\boldsymbol u}(\tau)- {\boldsymbol v}(\tau)\rangle
      _\sigma d\tau ~~~~~~~~~~~~~~~~~~~~~~~~~~$$
 $$  ~~+ \int_0^t \int_\Omega ({\boldsymbol u}(x,\tau)\cdot\nabla){\boldsymbol u}
(x,\tau)\cdot
      \nabla({\boldsymbol u}(x,\tau)-{\boldsymbol v}(x,\tau))dx d\tau 
     +\frac 12 |{\boldsymbol u}(t)-{\boldsymbol v}(t)|^2_{0,2}  \eqno{(1.2)}$$
  $$ ~~~~~~~~~~~~~~~\le \int_0^t ({\boldsymbol g}(\tau),{\boldsymbol u}(\tau)-{\boldsymbol v}
(\tau))_\sigma d\tau
     + \frac 12 |{\boldsymbol u}_0-{\boldsymbol v}(0)|^2_{0,2},
~~\forall t \in [0,T],
    ~\forall {\boldsymbol v} \in \boldsymbol{\cal K}(p).~~~~~~~~~~ $$
\end{description}
}

We note that $\bfu$ is defined for \emph{every} $t\in[0,T]$,
and according to the given $\bfu_0$,
even if we do not require it to be continuous in time:
our definition permits jumps in time, including 
the initial time $t=0$.
What we will prove, is that $\bfu$ is a limit of continuous
approximate solutions.

The main objective of 
this paper is to prove the following existence result for 
$P(p;{\boldsymbol g},{\boldsymbol u_0})$. 
\vspace{0.25cm}

\noindent
{\bf Theorem 1.1.} {\it Assume that (1.1) is satisfied and
 $${\boldsymbol u}_0 \in {\boldsymbol W}_\sigma(\Omega),~
      ~{\rm supp}({\boldsymbol u}_0)\subset \{x \in \Omega~|~p(x,0)>0\},~~
      |{\boldsymbol u}_0|\le p(\cdot,0)
   ~\text{ in~} \Omega.\eqno{(1.3)}$$ 
Then, there
is at least one solution ${\boldsymbol u}$
of $P(p;{\boldsymbol g},{\boldsymbol u_0})$. }\vspace{0.5cm}

We do not touch the uniqueness of solution problem, 
even if uniqueness holds
for constraints considered in [3--6, 13, 24, 25]. 
In our case it remains an open question, 
together with time continuity.
We state uniqueness for approximate solutions,
defined in Section~2 (see Proposition~2.1).

In the proof of Theorem~1.1 the main difficulty
arises from the nonlinear convection term 
$(\bfu\cdot \nabla)\bfu$. In our case, the class of test functions 
$\boldsymbol{\cal K}(p)$ is not a linear space. 
Therefore, the usual compactness methods, based on Sobolev embeddings,
cannot be directly applied. 
Our idea is to use local bounded variation estimate of $\bfu$ 
or its approximate solutions in the 
space $(0,T) \times {\boldsymbol W}_\sigma^*(\Omega)$; see Section~3. 
This is a completely new approach to 
parabolic variational inequalities of the Navier--Stokes type. 
The proof of Theorem~1.1 is given in Section~4.

\section{Approximate problems}\label{sec:appro}

In this section, we propose an approximation procedure to 
$P(p;{\boldsymbol g},{\boldsymbol u_0})$. We begin with a regular 
approximation of the obstacle
function $p(x,t)$. Choose a sequence $\{p_n\}$ of Lipschitz,
non-degenerate obstacle functions on~$\overline Q$ such that 
 $$ \left \{
  \begin{array}{l}
    0 < p_n(x,t) < \infty~\text{ on~} \overline Q,~\forall n \in \text{ N},
    \\[0.18cm]
     \forall \kappa \in (0,\infty),~ p_n \stackrel{n\to\infty}{\longrightarrow} p~
    \text{ uniformly~in~} \overline Q(p\le \kappa):=\{(x,t)\in 
     \overline Q~|~p(x,t)\le \kappa\},\\[0.18cm]
    \text{for~ any~ sufficiently~ large~}M >0, 
    \text{ there is an integer } n_M \in {\bf N}~\text{ such~that}\\
     ~~~M \le p_n \le p
~\text{ on~}\overline Q(p>M):=\{(x,t)\in \overline Q~|~p(x,t)>M\},
 ~\forall n\geq n_M.
    \end{array} \right.  \eqno{(2.1)} $$
{\bf Remark 2.1.} Given function $p(x,t)$ satisfying (1.1), we always
construct an approximate sequence $\{p_n\}$ satisfying (2.1). For instance,
the sequence $\{p_n\}$ consisting of regularizations of cut-off functions
  $$ \tilde p_n(x,t):=\left \{
      \begin{array}{ll}
     \displaystyle{   n,}&\displaystyle{~~~\text{ if~}p(x,t) >n,}\\[0.2cm]
      \displaystyle{  p(x,t),}&\displaystyle{~~~\text{ if~}\frac 1n 
\le p(x,t) \le n,}\\[0.2cm] 
     \displaystyle{   \frac 1n,}&\displaystyle{~~~\text{ if~}0\le p(x,t) 
   < \frac 1n,}
      \end{array} \right. $$
fulfills (2.1).\vspace{0.5cm}

Next, we formulate precisely the approximate problem, denoted by $P(p_n; {\boldsymbol g}, {\boldsymbol u}_{0n})$.
\vspace{0.25cm}

\noindent
{\bf Definition 2.1} {\it
Given a function $p_n$ satisfying (2.1) and an initial datum 
 $${\boldsymbol u}_{0n} \in {\boldsymbol V}_\sigma(\Omega),~~
  |{\boldsymbol u}_{0n}|\le p_n(\cdot,0)~\text{ a.e.~in~} \Omega, \eqno{(2.2)}$$
the problem $P(p_n; {\boldsymbol g}, {\boldsymbol u}_{0n})$ is
to find a function ${\boldsymbol u}_n=(u_n^{(1)}, u_n^{(2)},u_n^{(3)})$ which
satisfies the following (1) and (2):
\begin{description}
\item {(1)} ${\boldsymbol u}_n \in W^{1,2}(0,T;{\boldsymbol H}_\sigma(\Omega))
\cap C([0,T];{\boldsymbol V}_\sigma(\Omega))$ such that
  $$ |{\boldsymbol u}_n(x,t)| \le p_n(x,t)~\text{ a.e.~}x \in \Omega,~\forall
     t \in [0,T], 
     \eqno{(2.3)}$$
  $$ \begin{array}{l}
   \displaystyle{ ({\boldsymbol u}'_n(t), {\boldsymbol u}_n(t)-{\boldsymbol z})_\sigma
    +\nu \langle {\boldsymbol u}_n(t), {\boldsymbol u}_n(t)-{\boldsymbol z}
\rangle_\sigma~~~~~~~~~~~
    }\\[0.3cm]
  \displaystyle{+ \int_\Omega ({\boldsymbol u}_n(t)\cdot \nabla)
   {\boldsymbol u}_n(t)\cdot ({\boldsymbol u}_n(t)
        -{\boldsymbol z}) dx }
  \le ({\boldsymbol g}(t), {\boldsymbol u}_n(t)-{\boldsymbol z})_\sigma,
    \\[0.4cm]
  \forall {\boldsymbol z} \in {\boldsymbol V}_\sigma(\Omega)~\text{ with~}|{\boldsymbol z}|\le 
   p_n(\cdot,t)~\text{ a.e.~in~}\Omega,~~\text{ a.e.~}t \in [0,T].
     \end{array}  \eqno{(2.4)}$$ 
\item{(2)} ${\boldsymbol u}_n(0)={\boldsymbol u}_{0n}$ in ${\boldsymbol H}_\sigma(\Omega)$.
\end{description}
}

We are now applying the general theory on evolution inclusions
generated by time dependent subdifferentials to the solvability of 
$P(p_n; {\boldsymbol g},{\boldsymbol u}_{0n})$.
To this end, we introduce a time dependent convex function $\varphi^t_n,~
t \in [0,T],$ on
${\boldsymbol H}_\sigma(\Omega)$ given by:
 $$ \varphi^t_n({\boldsymbol z}):=\left \{
      \begin{array}{ll}
     \displaystyle{   \frac \nu 2 |{\boldsymbol z}|^2_{1,2} 
    +I_{K(p_n;t)}({\boldsymbol z}),}&\displaystyle{~~\forall {\boldsymbol z} \in 
{\boldsymbol V}_\sigma(\Omega),}\\[0.3cm]
      \displaystyle{ \infty,}&~~\text{ otherwise,} 
      \end{array} \right.  $$
where 
$$K(p_n;t):=\{{\boldsymbol z} \in {\boldsymbol V}_\sigma(\Omega)~
|~|{\boldsymbol z}| 
\le p_n(\cdot,t)~\text{ a.e.~in~}\Omega\},
$$
which is closed and convex in
${\boldsymbol V}_\sigma(\Omega)$, and $I_{K(p_n;t)}$ is its indicator
function on ${\boldsymbol V}_\sigma(\Omega)$, namely
 $$  I_{K(p_n;t)}({\boldsymbol z}):= \left\{
      \begin{array}{l}
          0,~~~~\text{ if ~}{\boldsymbol z} \in K(p_n;t),\\
          \infty,~~~\text{ otherwise}.
      \end{array}  \right.  $$
Clearly $\varphi^t_n$ is non-negative, proper, l.s.c. and convex on 
${\boldsymbol H}_\sigma(\Omega)$ and on ${\boldsymbol V}_\sigma(\Omega)$ 
for every 
$t \in [0,T]$. Also, we define a perturbation term 
${\boldsymbol G}(t,\cdot): K(p_n;t) \to {\boldsymbol H}_\sigma(\Omega)$ by the formula:
 $$ ({\boldsymbol G}(t,{\boldsymbol v}),{\boldsymbol z})_\sigma:= 
  \int_\Omega ({\boldsymbol v}\cdot \nabla){\boldsymbol v} \cdot 
   {\boldsymbol z}dx = \sum_{k,j=1}^3 
  \int_\Omega
     v^{(k)}(x)\frac {\partial v^{(j)}(x)}{\partial x_k}z^{(j)}(x) dx
  $$
for all ${\boldsymbol v}=(v^{(1)}, v^{(2)},v^{(3)})\in K(p_n;t)$ and 
${\boldsymbol z}=(z^{(1)}, z^{(2)},z^{(3)})\in {\boldsymbol H}_\sigma(\Omega)$.
\vspace{0.3cm}

\noindent
{\bf Lemma 2.1.} {\it Let ${\boldsymbol u}_n$ be a function in 
$W^{1,2}(0,T;{\boldsymbol H}_\sigma(\Omega))\cap C([0,T];{\boldsymbol V}_\sigma(\Omega))$ and let ${\boldsymbol u}_{0n} \in K(p_n;0)$. Then,
$P(p_n; {\boldsymbol g}, {\boldsymbol u}_{0n})$ is equivalent to the 
following Cauchy problem:
 $$ \left\{
    \begin{array}{l}
 \displaystyle{{\boldsymbol u}'_n(t)+\partial \varphi^t_n({\boldsymbol u}_n(t))
   +{\boldsymbol G}(t,{\boldsymbol u}_n(t))
       \ni {\boldsymbol g}(t)~~\text{ in~}{\boldsymbol H}_\sigma(\Omega),
    ~\text{ a.e.~}t \in [0,T],} \\[0.3cm]
      {\boldsymbol u}_n(0)={\boldsymbol u}_{0n},
     \end{array} \right. \eqno{(2.5)}$$
where $\partial \varphi^t_n(\cdot)$ is the subdifferential of 
$\varphi^t_n(\cdot)$ in ${\boldsymbol H}_\sigma(\Omega)$.}\vspace{0.3cm}

The equivalence required  in Lemma 2.1 is derived  immediately from the definition of
the subdifferential $\partial \varphi^t_n$, namely ${\boldsymbol v}^* \in \partial 
\varphi^t_n({\boldsymbol v})$ if and only if ${\boldsymbol v}^* \in 
{\boldsymbol H}_\sigma(\Omega),~
{\boldsymbol v}\in K(p_n;t)$ and
 $$ ({\boldsymbol v}^*, {\boldsymbol z}-{\boldsymbol v})_\sigma 
        +\nu \langle {\boldsymbol v},{\boldsymbol v}-{\boldsymbol z}\rangle_\sigma \le 0,
      ~~\forall z\in K(p_n,t). $$ 
For a detailed proof, see [13, 14].\vspace{0.3cm}

\noindent
{\bf Lemma 2.2.} {\it Let $p_n$ satisfy (2.1).
There is a positive constant $C_n$, depending 
on $n$, such that for every $s,t \in [0,T]$ and every 
${\boldsymbol z} \in K(p_n;s)$ there is $\tilde{\boldsymbol z} \in K(p_n;t)$ satisfying 
 $$|\tilde {\boldsymbol z}-{\boldsymbol z}|_{0,2} 
     \le C_n|p_n(t)-p_n(s)|_{C(\overline \Omega)},~~\varphi^t_n(\tilde{\boldsymbol z})
      \le \varphi^s_n({\boldsymbol z}). \eqno{(2.6)}$$ 
}

\noindent
{\bf Proof.} 
Let  $\mu_n = \min_{(x,t)\in\bar{Q}}p_n(x,t)$; 
we note that $\mu_n>0$ by (2.1). 
Denote also by $L_n$ the Lipschitz constant of $p_n$ and 
choose a partition of $[0,T]$,
$0=t_0<t_1<t_2<\cdots <t_N=T$,  so that 
$t_{k}-t_{k-1}\le{\mu_n}/{L_n}$ for $k=1,\ldots,N$.
Then
  $$ |p_n(t)-p_n(s)| <\mu_n,
     \quad\text{ for } t,s\in[t_k-t_{k-1}],~~k=1,2,\cdots,N. 
  $$
Now, suppose that $s,t \in [t_{k-1},t_k]$. 
Given ${\bfz} \in K(p_n;s)$, we consider the function
$$\tilde {\bfz}(x):=\left(1-\frac 1\mu_n
     |p_n(s)-p_n(t)|
        _{C(\overline \Omega)} \right){\bfz}(x). $$
Observe that
 \begin{eqnarray*}
    |\tilde{\bfz}(x)|
      &=& 
      \left(1-\frac 1\mu_n |p_n(s)-p_n(t)|_{C(\overline \Omega)} \right)
      |{\bfz}(x)|\\
     &\le & 
     \left(1-\frac 1\mu_n |p_n(s)-p_n(t)|_{C(\overline \Omega)} \right)
     p_n(x,s)\\
     &\le & p_n(x,s)-|p_n(s)-p_n(t)|_{C(\overline \Omega)}\\
     &\le & p_n(x,s)-|p_n(x,s)-p_n(x,t)|
     \  \le \  p_n(x,t).
 \end{eqnarray*} 
Since div $\tilde{\bfz}=0$, the above inequality implies that
$\tilde{\bfz} \in K(p_n;t)$. Clearly, the second inequality of~(2.6) is satisfied.
Moreover,
   $$ |\tilde{\bfz}- {\bfz}|_{0,2} =
   \frac{1}{\mu_n}\,|p_n(s)-p_n(t)|_{C(\overline \Omega)}|{\bfz}|_{0,2}\\
   \le C_n' \,|p_n(t)-p_n(s)|_{C(\overline \Omega)}
   $$
where 
$$
C_n' = \frac{\max\limits_{(x,t)\in\bar{Q}}p_n(x,t)}{\mu_n}\,|\Omega|^\frac{1}{2} ,
$$
which is finite by~(2.1).

In the general case of $s,t \in [0,T]$ and 
${\bfz}\in K(p_n;s)$,
by repeating the above procedures at most $N$-times, we can construct
$\tilde {\bfz} \in K(p_n;t)$ satisfying both required inequalities, the first one 
with the constant
$C_n:= NC_n'$.  \hfill $\Box$ \vspace{0.5cm}

We prepare a lemma which we shall need in the convergence of approximate 
problems.\vspace{0.5cm}

\noindent
{\bf Lemma 2.3.} {\it Let $p_n$ satisfy (2.1).
Let ${\boldsymbol v}$ be any function in
$C([0,T]; {\boldsymbol W}_\sigma(\Omega))$ such that 
$$ {\rm supp}({\boldsymbol v})\subset \{(x,t) \in \Omega\times [0,T]
  ~|~p(x,t)\geq \delta \},~~|{\boldsymbol v}| \le p~\text{ on}~Q
$$
for a positive number $\delta$, and put
 $$ \delta_n := \frac 1\delta \max_{\text{ supp}({\boldsymbol v})}
      |p\land M-p_n\land M|,~~\forall n \in {\bf N}, $$
where $M=\delta+\sup_{(x,t) \in Q} |{\boldsymbol v}(x,t)|$, 
$p \land M= \min\{p,M\}$ and $p_n \land M= \min\{p_n,M\}$.
Then $\delta_n\to 0 $ and for any $n \in {\bf N}$, we have 
 $$(1-\delta_n)^+{\boldsymbol v}(t) \in K(p_n; t),~~\forall t \in [0,T].
   \eqno{(2.7)}$$ }

\noindent
{\bf Proof.} It follows easily from (2.1) that $\delta_n \to 0$ as
$n \to \infty$. In addition, in case $\delta_n <1$, 
  \begin{eqnarray*}
    (1-\delta_n)^+|{\boldsymbol v}(x,t)| 
   &\le & 
      (1-\delta_n)|{\boldsymbol v}(x,t)|
    \ \le \  
      (1-\delta_n)\,p(x,t)\land M\\
    & = & p(x,t)\land M-\frac {p(x,t)\land M}\delta |p(x,t)\land M
      -p_n(x,t)\land M|\\
    & \le&  p(x,t)\land M -(p(x,t)\land M-p_n(x,t)\land M) \le p_n(x,t)
  \end{eqnarray*}
for all $(x,t) \in {\rm supp}({\boldsymbol v})$. This shows (2.7).
 \hfill $\Box$
\vspace{0.5cm}

\noindent
{\bf Proposition 2.1.} {\it Let $p_n$ satisfy (2.1).
Let ${\boldsymbol u}_0$ be any element in
${\boldsymbol W}_\sigma(\Omega)$ for which (1.3) holds; hence 
$${\rm supp}(\bfu_0) \subset \{x \in \Omega~|~p(x,0)\geq \hat\delta\}
$$ for a certain
constant $\hat\delta>0$. Put
$$ {\boldsymbol u}_{0n} :=(1-\hat\delta_n)^+{\boldsymbol u}_0~~{\it with~}
   \hat\delta_n:=\frac 1{\hat\delta} \max_{{\rm supp}(\bfu_0)}
   |p\land \hat M-p_n\land \hat M|,
  ~~\forall n \in {\bf N},$$
where $\hat M=\hat\delta+|\bfu_0|_{C(\overline \Omega)}$. 
Then problem $P(p_n;{\boldsymbol g}, {\boldsymbol u}_{0n})$ admits one and 
only one solution 
${\boldsymbol u}_n$ in $W^{1,2}(0,T;{\boldsymbol H}_\sigma(\Omega))
\cap C([0,T];{\boldsymbol V}_\sigma(\Omega))$, which is also the unique solution of
(2.5). Moreover,
${\boldsymbol u}_n$ satisfies the estimate
$$ |{\boldsymbol u}_n(t)|^2_{0,2} +\nu \int_0^t |{\boldsymbol u}_n|^2_{1,2}
  d\tau
   \le |{\boldsymbol u}_0|^2_{0,2}+ \frac {L_P^2}{\nu}\int_0^T
   |{\boldsymbol g}|^2_{0,2} d\tau =: M_0,
  ~\forall t \in [0,T], \eqno{(2.8)}$$ 
where $L_P$ is the Poincar\'e constant, i.e.
 $$|{\boldsymbol z}|_{0,2}
  \le L_P|{\boldsymbol z}|_{1,2},~~ \forall {\boldsymbol z} \in 
    {\boldsymbol V}_\sigma.$$
}

\noindent
{\bf Proof.} We show in a similar way to that of Lemma 2.3 that
${\boldsymbol u}_{0n}$ satisfies (2.2) for all~$n$.
The time dependence (2.6) 
of the mapping $t \mapsto \varphi^t_n$ 
is a sufficient condition for the Cauchy 
problem (2.5) without perturbation ${\boldsymbol G}$ 
to have one and only one 
solution ${\boldsymbol u}$ (see [17, 22, 28]).
Furthermore, as to the perturbation term ${\boldsymbol G}$, we have 
 $$ |({\boldsymbol G}(t,{\boldsymbol v})-{\boldsymbol G}(t,\bar{\boldsymbol v}), {\boldsymbol v}-\bar {\boldsymbol v})_\sigma|
    \le \varepsilon |{\boldsymbol v}-\bar{\boldsymbol v}|^2_{1,2}
      +C_\varepsilon |{\boldsymbol v}-\bar{\boldsymbol v}|^2_{0,2},
    ~~\forall {\boldsymbol v},~\bar{\boldsymbol v} \in K(p_n;t),\eqno{(2.9)}$$
where $\varepsilon$ is any positive constant and $C_\varepsilon$ is a positive 
constant depending only on~$\varepsilon$ and~$n$.
Indeed, 
\begin{multline*}
({\boldsymbol G}(t,{\boldsymbol v})-{\boldsymbol G}(t,\bar{\boldsymbol v}), {\boldsymbol v}-\bar {\boldsymbol v})_\sigma
=
\sum_{k,j=1}^3
\int_\Omega \left( v^{(k)}
\frac{\partial v^{(j)}}{\partial x_k}-\bar{v}^{(k)} 
\frac{\partial \bar{v}^{(j)}}{\partial x_k}\right)
(v^{(j)}-\bar{v}^{(j)})dx
\\
=
\sum_{k,j=1}^3
\int_\Omega (v^{(k)} -\bar{v}^{(k)} )
\frac{\partial v^{(j)}}{\partial x_k}(v^{(j)}-\bar{v}^{(j)})dx
+
\sum_{k,j=1}^3
\int_\Omega \bar{v}^{(k)} 
\frac{\partial (v^{(j)}-\bar{v}^{(j)})}{\partial x_k}(v^{(j)}-\bar{v}^{(j)})dx,
\end{multline*}
and from the fact that $\text{div}\,\bfv = \text{div}\,\bar{\bfv}=0$ we infer
that the second sum is equal to~$0$, while the first is estimated by 
$ 9(\max_\Omega |\bfv|)\,|{\boldsymbol v}-\bar{\boldsymbol v}|_{1,2} 
|{\boldsymbol v}-\bar{\boldsymbol v}|_{0,2}$, so that (2.9) follows. 
Therefore, according to the perturbation result of~[26],
$P(p_n; {\boldsymbol g}, {\boldsymbol u}_{0n})$ has a unique solution 
${\boldsymbol u}_n$. Also, by taking ${\boldsymbol z}=0$ in (2.4) and 
integrating in time over $[0,t]$ we get
$$ \frac 12 |{\boldsymbol u}_n(t)|^2_{0,2} 
  +\nu \int_0^t |{\boldsymbol u}_n|^2_{1,2} d\tau
   \le \frac 12 |{\boldsymbol u}_{0n}|^2_{0,2} 
  + \int_0^t ({\boldsymbol g},{\boldsymbol u}_n)_\sigma d\tau. \eqno{(2.10)}$$
Noting that
  $$ \int_0^t |({\boldsymbol g},{\boldsymbol u}_n)_\sigma| d\tau
   \le \frac \nu 2 \int_0^t |{\boldsymbol u}_n|^2_{1,2}d\tau
    + \frac {L^2_P}{2\nu} \int_0^t |{\boldsymbol g}|^2_{0,2} d\tau,$$
we obtain (2.8) from (2.10), since 
$|{\boldsymbol u}_{0n}|_{0,2} \le |{\boldsymbol u}_0|_{0,2}$. \hfill $\Box$
 \vspace{0.5cm}

It follows from the energy estimate (2.8) that there exists 
a subsequence $\{{\boldsymbol u}_{n_k}\}$ and a function 
$\bfu\in L^2(0,T; {\boldsymbol V}_\sigma(\Omega))$ such that 
${\boldsymbol u}_{n_k}$ weakly converges to ${\boldsymbol u}$ in
 $L^2(0,T;{\boldsymbol V}_\sigma(\Omega))$ as well as weakly$^*$ in
$L^\infty(0,T; {\boldsymbol H}_\sigma(\Omega))$ (as $k\to \infty$). In the 
sequel, for simplicity of notation
we write $\{{\boldsymbol u}_n\}$ again for $\{{\boldsymbol u}_{n_k}\}$, namely 
 $$ {\boldsymbol u}_n \to {\boldsymbol u}
 ~\text{ weakly~in~}L^2(0,T;{\boldsymbol V}_\sigma(\Omega))~
    \text{ and weakly}^*~\text{ in }L^\infty(0,T; {\boldsymbol H}_\sigma(\Omega)) \text{ as } n\to\infty.\eqno{(2.11)} $$
We will refer to~(2.11) in most statements and proofs which follow.
We underline however that we have here only a subsequence 
of the sequence $\bfu_n$ constructed in Proposition~2.1.

\section{Local uniform estimate of the total variation of ${\boldsymbol u}_n$}\label{sec:bv}

In this section we use the notation:
\begin{eqnarray*}
&& \hat Q:=\Omega\times [0,T],\\
&& \hat Q(p> \kappa):=\{(x,t)\in \hat Q~|~p(x,t)>\kappa\},~~0<\kappa<\infty,\\
&& \hat Q(p=\infty):=\{(x,t)\in \hat Q~|~p(x,t)=\infty\},\\
&&\hat Q(p=0):=\{(x,t)\in \hat Q~|~p(x,t)=0\};
\end{eqnarray*}
by (1.1), $\hat Q(p> \kappa)$ is relatively open in $\hat{Q}$, 
and $\hat Q(p= \infty)$ and $\hat Q(p=0)$ 
are relatively compact in $\hat{Q}$.
We will also use in this section the spaces $\bfW_\sigma(\Omega'), \, 
\bfV_\sigma(\Omega'), \, \bfH_\sigma(\Omega')$, built 
on any open set $\Omega'\subset\Omega$, and use the same notation 
as in Section~\ref{sec:intro} for the norms 
without indicating $\Omega'$ explicitly therein.

We shall use the continuous embeddings:
 $$ {\boldsymbol W}_\sigma(\Omega') \hookrightarrow C(\overline {\Omega'})^3,~~
    {\boldsymbol W}_\sigma(\Omega') \hookrightarrow 
     {\boldsymbol V}_\sigma(\Omega') \hookrightarrow L^4(\Omega')^3
     \hookrightarrow L^2(\Omega')^3,$$
with inequalities
\begin{equation}\label{emb3}
\begin{array}{c}
|{\boldsymbol  f}|_{C(\overline{\Omega'})^3} \le L_0 \, 
|{\boldsymbol f}|_{1,4},  
\quad
   |{\boldsymbol  f}|_{1,2} \le 
   L_2 \, |{\boldsymbol f}|_{1,4},
\quad 
\forall {\boldsymbol f}\in {\boldsymbol  W}_\sigma(\Omega'), \\[0.3cm]
 |{\boldsymbol  f}|_{0,2} \le L_1 |{\boldsymbol  f}|_{1,2},~~
|{\boldsymbol f}|_{0,4}:=
|{\boldsymbol f}|_{L^4(\Omega')^3} \le  L_3 \, |{\boldsymbol f}|_{1,2},
\quad  \forall {\boldsymbol f}\in {\boldsymbol  V}_\sigma(\Omega'),
\end{array}
\end{equation}
where $L_0,~L_1:=L_P,~L_2$ and $L_3$ are positive constants,
which are derived from the Sobolev inequalities 
as well as the Poincar\'e inequality, 
and independent of $\Omega' \subset \Omega$ (cf. [15; Ch.7]).\
\vspace{0.3cm}

\noindent
{\bf Lemma 3.1.} {\it 
Let $p$ satisfy (1.1), $\kappa>0$, $\Omega'$ be an open set in $\Omega$, and
$\hat Q':=\Omega'\times [T_1,T_1']$. Assume that 
$\tilde\Omega'\times [T_1,T_1'] \subset \hat Q(p>\kappa)$,
where $\tilde \Omega'$ is the relative closure of 
$\Omega'$ in $\Omega$. Then, for  $n$ large enough,  ${\boldsymbol u}_n$
defined in Proposition~2.1
is of bounded variation as a function
from $[T_1,T_1']$ into ${\boldsymbol W}^*_\sigma(\Omega')$. 
Its total variation is uniformly bounded with respect to $n$ 
and the bound depends only on $\kappa$.}\vspace{0.3cm}

\noindent
{\bf Proof.} 
By (2.1), we fix $N = N(\kappa)$ large enough to have
\begin{equation}\label{eq:700}
\frac{\kappa}{2}\le p_n(x,t), \quad  \forall n>N, \ \forall (x,t)\in \hat Q(p>\kappa).
\end{equation}
Take now ${\boldsymbol z}\in C([T_1, T'_1]; {\boldsymbol W}_\sigma(\Omega'))$ 
with 
\begin{equation}\label{eq:z}
 {\boldsymbol z}(T_1)={\boldsymbol z}(T'_1)=0,~~
\text{supp }(\bfz)\subset \Omega'\times [T_1,T'_1], ~~ 
\|{\boldsymbol z}\|_{C([T_1,T_1'];{\boldsymbol W}_\sigma(\Omega'))} \le 
\frac{\kappa}{2L_0},
\end{equation}
where $L_0$ is defined by~\eqref{emb3}. Then by~\eqref{emb3}, 
\eqref{eq:700} and~\eqref{eq:z},
$$
|{\boldsymbol z}(x,t)|\le \|{\boldsymbol z}\|_{C([T_1,T'_1];{C(\overline
{\Omega'})}^3)}\le 
\frac{\kappa}{2} \le p_n(x,t), 
\quad \forall (x,t)\in\hat{Q}(p>\kappa), \quad \forall n>N.
$$
So, $\bfz(\cdot,t)\in K(p_n;t)$ for $n> N$, i.e.~$\bfz(t)$ is a proper test function in~(2.3), which writes
\begin{multline*}
  ({\boldsymbol u}'_n(t), {\boldsymbol u}_n(t)-{\boldsymbol z}(t))_\sigma
    +\nu \langle {\boldsymbol u}_n(t), {\boldsymbol u}_n(t)-{\boldsymbol z(t)}\rangle_\sigma 
+ (\bfG(t,{\boldsymbol u}_n(t)), {\boldsymbol u}_n(t)-{\boldsymbol z}(t))_\sigma \\
  \le ({\boldsymbol g}(t), {\boldsymbol u}_n(t)-{\boldsymbol z}(t))_\sigma,\quad\text{ for a.~e.~} t\in[T_1,T_1'].
 \end{multline*}
As $ (\bfG(t,{\boldsymbol u}_n(t)), {\boldsymbol u}_n(t))_\sigma = 0$, we have for a.e.~$t\in[T_1,T'_1]$, all $n>N$
and ${\boldsymbol z}$ satisfying \eqref{eq:z}:
\begin{multline*}
  -({\boldsymbol u}'_n(t), {\boldsymbol z}(t))_\sigma 
   +\frac{1}{2}\frac{d}{dt} |{\boldsymbol u}_n(t)|^2_{0,2}
 +\nu |{\boldsymbol u}_n(t)|_{1,2}^2  \\
\le
    \nu |{\boldsymbol u}_n(t)|_{1,2} |{\boldsymbol z}(t)|_{1,2} 
+
 (\bfG(t,{\boldsymbol u}_n(t)), {\boldsymbol z}(t))_\sigma + 
  ({\boldsymbol g}(t), {\boldsymbol u}_n(t)-{\boldsymbol z}(t))_\sigma.
 \end{multline*}
When integrated in time, with the Young, Schwarz and Poincar\'e inequalities, this implies:
$$\begin{array}{l}
\displaystyle
  -\int_{T_1}^{T_1'}({\boldsymbol u}'_n(t), {\boldsymbol z}(t))_\sigma \,dt 
    +\frac{1}{2} |{\boldsymbol u}_n(T_1')|^2_{0,2}
 +\frac{\nu}{2} \int_{T_1}^{T_1'} |{\boldsymbol u}_n(t)|_{1,2}^2\,dt \\[0.15cm]\displaystyle
~~~~~~~~~~~\le 
 \int_{T_1}^{T_1'}(\bfG(t,{\boldsymbol u}_n(t)), {\boldsymbol z}(t))_\sigma\, dt 
 + \nu \int_{T_1}^{T_1'} |{\boldsymbol u}_n(t)|_{1,2} |{\boldsymbol z}(t)|_{1,2}\, dt\\[0.15cm]\displaystyle
 ~~~~~~~~~~~~~~~~~~~~~~
  + \int_{T_1}^{T'_1}|{\boldsymbol g}(t)|_{0,2} 
  |{\boldsymbol z}(t)|_{0,2} \,dt
+\frac{1}{2} |{\boldsymbol u}_n(T_1)|^2_{0,2} 
+ \frac{L_1^2}{2\nu}\int_{T_1}^{T_1'} |{\boldsymbol g}(t)|_{0,2}^2\,dt . 
\end{array}\eqno{(3.4)}
$$
Note that
$$
   \nu \int_{T_1}^{T'_1} |{\boldsymbol u}_n|_{1,2}
  |{\boldsymbol z}|_{1,2}dt
 \le \nu 
   \|{\boldsymbol u}_n\|_{L^2(0,T;{\boldsymbol V}_\sigma(\Omega))}
    \|{\boldsymbol z}\|_{L^2(T_1,T'_1;{\boldsymbol V}_\sigma(\Omega'))},
 \eqno{(3.5)}$$
$$
   \int_{T_1}^{T'_1} |{\boldsymbol g}|_{0,2}
  |{\boldsymbol z}|_{0,2}dt
   \le 
    \|{\boldsymbol g}\|_{L^2(0,T;{\boldsymbol H}_\sigma(\Omega))}
 \|{\boldsymbol z}\|_{L^2(T_1,T'_1;{\boldsymbol H}_\sigma(\Omega'))}~~~~$$ 
 $$~~~~~~~~~~~~~~~~~~~~~  \le L_1 \|{\boldsymbol g}\|_{L^2(0,T;{\boldsymbol H}
    _\sigma(\Omega))}
    \|{\boldsymbol z}\|_{L^2(T_1,T'_1;{\boldsymbol V}_\sigma(\Omega'))}.
  \eqno{(3.6)}
 $$
Besides, as div ${\boldsymbol u}_n = $ div ${\boldsymbol z} = 0$, we have
\begin{multline*}
(\bfG(t,{\boldsymbol u}_n(t)), {\boldsymbol z}(t))_\sigma  
\le 
\left|
\sum\limits_{k,j = 1}^{3}
\int_\Omega u_n^{(k)}(x,t)\frac{\partial u_n^{(j)}(x,t)}{\partial x_k}z^{(j)}(x,t)\:dx
\right|
\\
\le 
\left|\sum\limits_{k,j = 1}^{3}
\int_\Omega u_n^{(k)}(x,t)\frac{\partial z_n^{(j)}(x,t)}{\partial x_k}u_n^{(j)}(x,t)\:dx
\right| 
\ 
\le 9 |{\boldsymbol u}_n(t)|_{0,2} | {\boldsymbol u}_n(t)|_{0,4}
 |{\boldsymbol z}(t)|_{1,4},
\end{multline*}
so that
$$\begin{array}{lrl}
& \displaystyle{\int_{T_1}^{T_1'}(\bfG(t,{\boldsymbol u}_n(t)), 
  {\boldsymbol z}(t))_\sigma dt}
\le &
9  L_3 \, \|{\boldsymbol u}_n\|_{L^\infty(0,T;{\boldsymbol H}_\sigma(\Omega))}
\| {\boldsymbol u}_n\|_{L^2(0,T;{\boldsymbol V}_\sigma(\Omega))} 
 \|{\boldsymbol z}\|_{L^2(T_1,T'_1;{\boldsymbol W}_\sigma(\Omega'))}
\end{array} \eqno{(3.7)} $$
Put (3.5)--(3.7) into (3.4) and neglect the positive terms 
at the left hand side. Then from (2.8) we obtain
for all ${\boldsymbol z}$ satisfying \eqref{eq:z} and all $n>N$ 
$$
\begin{array}{ll}
   \displaystyle{-\int_{T_1}^{T_1'}({\boldsymbol u}'_n(t), 
 {\boldsymbol z}(t))_\sigma dt} 
&\displaystyle{\le 
 M_0+ M_1 \|{\boldsymbol z}\|_{L^2(T_1,T'_1;{\boldsymbol V}_\sigma(\Omega'))}
+M_2\|{\boldsymbol z}\|_{L^2(T_1,T'_1;{\boldsymbol W}_\sigma(\Omega'))}}
\\[0.3cm]
&\displaystyle{\le M_0+ M_3\|{\boldsymbol z}\|_{L^2(T_1,T'_1;
  {\boldsymbol W}_\sigma(\Omega'))}}\\[0.2cm]
&\displaystyle{\le M_0+ M_3T^{\frac 12} \|{\boldsymbol z}\|
  _{C([T_1,T'_1];\bfW_\sigma(\Omega'))},}
\end{array} \eqno{(3.8)}$$
where 
 $$
 M_1:=
 (\nu M_0)^{\frac 12}
 +
 L_1\|{\boldsymbol g}\|_{L^2(0,T;{\boldsymbol H}_\sigma(\Omega))},~ M_2:=9L_3 \nu^{-\frac 12}M_0, ~
     M_3:=M_1 L_2^{\frac 12}+M_2.$$
Since  $-\bfz$ is also a possible test function, 
we actually have that for all 
${\boldsymbol z}$ satisfying \eqref{eq:z} and $n>N$
\begin{eqnarray*}
   \left|
   \int_{T_1}^{T_1'}
   \left({\boldsymbol u}'_n(t),{\boldsymbol z}(t)\right)_\sigma dt 
   \right| 
   & \le & 
   M_0+M_3 T^{\frac 12} 
   \|{\bfz}\|_{C([T_1,T_1'];{\bfW}_\sigma(\Omega'))}.
\end{eqnarray*}
Take finally any 
$\tilde{\bfz}
\in 
C^1_0(T_1,T_1';{\boldsymbol W}_\sigma(\Omega'))$,  
put 
$$
{\boldsymbol z} = \frac{\tilde{\bfz}}{\|{\tilde{\bfz}}\|_{L^\infty(T_1,T_1';{\boldsymbol W}_\sigma(\Omega'))}}\cdot\frac{\kappa}{2 L_0}
$$
into the above inequality, and obtain that
$$
\left|\int_{T_1}^{T_1'}({\boldsymbol u}_n(t), 
{\tilde{\bfz}}'(t))_\sigma dt 
\right| 
=  \left|\int_{T_1}^{T_1'}({\boldsymbol u}'_n(t), {\tilde{\bfz}}(t))_\sigma dt \right| 
\le 
M_\kappa \,
\|{\tilde{\bfz}}\|_{L^\infty(T_1,T_1';{\bfW}_\sigma(\Omega'))}
$$
with
 $$ M_\kappa := \frac{2 L_0M_0}\kappa+M_3T^{\frac 12}, $$
for all ${\tilde{\bfz}}\in C^1_0(T_1,T_1';{\bfW}_\sigma(\Omega'))$ and all $n>N$. 
By a classical result on the relationship between weak derivatives and total 
variation, see e.g.~[7; Prop.~A.5], 
this implies that $\bfu_n \in BV(T_1,T'_1; \bfW_\sigma^*(\Omega'))$ 
and its total variation is bounded by $M_\kappa$.
\qed
\vspace{0.5cm}

\noindent
{\bf Lemma 3.2. }{\it  
Let $p$ satisfy (1.1) and
$\kappa>0$. Let $\Omega'$ be an open set 
in $\Omega$ such that
$\tilde \Omega'\times [T_1,T_1'] \subset \hat Q(p>\kappa)$,
where $\tilde \Omega'$ is the relative closure of 
$\Omega'$ in $\Omega$. Then, there exists a function $\bfu_{\Omega'}: 
[T_1,T'_1] \to \bfH_\sigma(\Omega')$, with 
$\sup_{T_1\le t \le T'_1} |{\bfu}_{\Omega'}(t)|_{\bfH_\sigma(\Omega')}\le M_0$ 
for the same constant
$M_0$ as in (2.8), such that, for $\bfu_n$ defined by~(2.11), 
$$\int_{\Omega'}{\bfu}_n(t)\cdot {\boldsymbol {\mathcal \xi}} dx  \to 
\int_{\Omega'}{\bfu}_{\Omega'}(t)
\cdot {\boldsymbol {\mathcal \xi}} dx,~{\it as~}n \to \infty, ~\forall 
  {\boldsymbol {\mathcal \xi}} \in {\bfH}_\sigma(\Omega'),~\forall t 
\in [T_1,T'_1]. $$ 
Moreover, $|\bfu_{\Omega'}(x,t)| \le p(x,t)$ for a.e.~$x \in \Omega'$ and every
$t \in [T_1,T'_1]$ and $\bfu_{\Omega'}=\bfu$ a.e.~on~$\Omega'\times [T_1,T'_1]$.}\vspace{0.5cm}

\noindent
{\bf Proof.} The space ${\bfW}_\sigma(\Omega')$ is separable. Let 
$X_0$ be its countable dense subset. At a first time, fix  any 
$\boldsymbol{\mathcal \xi}\in X_0$. We consider the 
sequence of real functions $f_n:[T_1, T_1']\to {\bf R}$ defined by
 $$f_n(t) = ({\bfu}_n(t), \boldsymbol{\mathcal \xi})_\sigma 
~\left(= \int_{\Omega'}
{\bfu}_n(x,t)\cdot \boldsymbol{\mathcal \xi}(x) dx\right).$$ 
Then, by Lemma~3.1, $f_n$ is uniformly bounded in 
$W^{1,1}(T_1, T_1')$, so is its total variation. Therefore,
it follows from the Helly selection theorem 
(see e.g.~[12; Section 5.2.3]) 
that there exists a subsequence $f_{n_k}$ of $\{f_n\}$ and a function 
$f\in BV(T_1,T_1')$
such that $f_{n_k}\to f$ pointwise on $[T_1,T'_1]$ and in $L^1(T_1,T_1')$.

However, the limit function $f$ and the subsequence $f_{n_k}$ depend 
also on $\boldsymbol{\mathcal \xi}$, that is, $n_k = n_k(\boldsymbol{\mathcal \xi})$ and  $f(t) = f(t;\boldsymbol{\mathcal \xi})$. 
But, since these are countable, by a diagonal argument we choose a 
subsequence, denoted again by $\{n_k\}$, such that  
$({\bfu}_{n_k}(t), \boldsymbol{\mathcal \xi})_\sigma$ converges to 
$f(t;\boldsymbol{\mathcal \xi})$
for all $\boldsymbol{\mathcal \xi}\in X_0$.
Furthermore, by density, this convergence holds for all 
$\tilde{\boldsymbol{\mathcal \xi}}\in{\bfW}_\sigma(\Omega')$. Indeed,
given  any $\eps>0$ and any 
$\tilde{\boldsymbol{\mathcal \xi}}\in{\bfW}_\sigma(\Omega')$, 
there exists $\boldsymbol{\mathcal \xi}\in X_0$
such that $|\boldsymbol{\mathcal \xi}-\tilde{\boldsymbol{\mathcal \xi}}|_{1,4} < \eps$, so that with 
$M =\limsup |\bfu_{n}|_{-1,\frac 43}~(<\infty)$ (cf. (2.8))
$$
\left|({\boldsymbol u}_{n_k}(t), \boldsymbol{\mathcal \xi})_\sigma - 
({\bfu}_{n_k}(t), \tilde{\boldsymbol{\mathcal \xi}})_\sigma\right|
\le
 |\bfu_{n_k}|_{-1,\frac 43}| \boldsymbol{\mathcal \xi} - 
\tilde{\boldsymbol{\mathcal \xi}}|_{1,4}
\le M\eps.
$$
This shows that $f(t,\tilde{\boldsymbol{\mathcal \xi}})
=\lim_{k\to \infty}({\boldsymbol u}_{n_k}(t), 
\tilde{\boldsymbol{\mathcal \xi}})_\sigma$ exists, i.e.\  we can extend $f(t;\cdot)$ to 
${\boldsymbol W}_\sigma(\Omega')$. 

Note that $f(t,\cdot)$ is linear as limit of linear functions. 
Thus, by the Riesz theorem, 
$$
\exists \tilde{\bfu}:[T_1,T_1']\to {\bfW}_\sigma^*(\Omega') 
\text{ such that }
f(t;{\boldsymbol{\mathcal \xi}}) = (\tilde{\bfu}(t), 
{\boldsymbol {\mathcal \xi}})
_\sigma, \quad\forall {\boldsymbol{\mathcal \xi}} \in {\bfW}_\sigma(\Omega').
$$
Since $|\bfu_{n_k}(t)|_{0,2}\le M_0$ for all $t \in [0,T]$, we have 
$\tilde{\bfu}(t) \in {\bfH}_\sigma(\Omega')$ for
all $t \in [T_1,T'_1]$ and
$\sup_{T_1\le t \le T'_1} |\tilde \bfu(t)|_{\bfH_\sigma(\Omega')}\le M_0$
and $\bfu_{n_k}(t) \to \tilde{\bfu}(t)$ weakly in 
$\bfH_\sigma(\Omega')$ as $k \to \infty$. 

We show now that $|\tilde{\bfu}|$ is bounded a.e.~by $p$. Take any $\eps>0$ and an integer 
$k(\eps)$ large enough to have $p_{n_k}(x,t) \le p(x,t)+\eps$ for all $x \in Q$ and
all $k \geq k(\eps)$, cf.~(2.1). Then $|\bfu_{n_k}(x,t)|\le p_{n_k}(x,t)
\le p(x,t)+\eps$ for a.e.~$(x,t) \in \Omega$. We note that the set 
$$F_{\eps}(t) = \{\bfz\in\bfH_\sigma(\Omega')~|~ |\bfz(x)|\le p(x,t)+\eps\text{ on } \Omega'\}
$$
 is convex and closed in~$\bfH_\sigma(\Omega')$. 
It follows from the Mazur lemma
(cf. [29; Th.2, Ch.V]) that the weak limit $\tilde{\bfu}(t)$ of 
$\bfu_{n_k}(t)$ in ${\bfH}_\sigma(\Omega')$ belongs to $F_\eps(t)$. 
By arbitrariness of 
$\eps>0$, we have
$|\tilde{\bfu}(x,t)| \le p(x,t)$ for a.e.~$x\in \Omega'$. 

We finally show that  $\tilde{\bfu} = \bfu$ a.e.~in $\Omega'\times [T_1,T'_1]$.
Take $\{E_i\}_{i=1}^N$ a partition of $[T_1,T_1']$ (family of pairwise disjoint measurable sets covering the interval) 
 and let $\boldsymbol{\mathcal \zeta}$ be a function of the form
$$\boldsymbol{\mathcal \zeta}(t)= \sum_{i=1}^N \chi_{E_i}(t)\boldsymbol{\mathcal \xi}_i \text{ with }  \boldsymbol{\mathcal \xi}_i\in\bfW_\sigma(\Omega').
$$
Then, by definition of $\tilde{\bfu}$ and since the sum is finite,
$$
\int_{T_1}^{T_1'} (\bfu_{n_k}(t),\boldsymbol{\mathcal \zeta}(t))_\sigma dt \to \int_{T_1}^{T_1'}
 (\tilde{\bfu}(t),\boldsymbol{\mathcal \zeta}(t))_\sigma dt.
$$
On the other hand, since 
$\bfu_n\to \bfu$ weakly in $L^2(0,T, \bfV_\sigma(\Omega'))$, cf.~(2.11),
and as $\bfW_\sigma(\Omega')\hookrightarrow\bfV_\sigma(\Omega')$, 
$$
\int_{T_1}^{T_1'} (\bfu_n(t),\boldsymbol{\mathcal \zeta}(t))_\sigma dt \to \int_{T_1}^{T_1'} ({\bfu}(t),\boldsymbol{\mathcal \zeta}(t))_\sigma dt.
$$
Consequently, 
by density of 
$\bfW_\sigma(\Omega')$ in $\bfH_\sigma(\Omega')$, 
$\tilde \bfu= \bfu$ a.e.~on 
$\Omega'\times [T_1,T'_1]$.
By uniqueness of the limit, all convergences stated above are valid for all the sequences, without extracting any subsequence. We obtain the statement of the lemma, where 
$\bfu_{\Omega'}:=\tilde {\bfu}$ is the required function.
\qed
\vspace{0.5cm}

\noindent
{\bf Corollary 3.1.} {\it Assume (1.1). Let $\Omega(t,0):=\{x \in \Omega~|~p(x,t)>0\}$
for each $t \in [0,T]$. Then there exists a function 
$\bar {\bfu}: [0,T] \to \bfH_\sigma(\Omega)$,
with $\sup_{t \in [0,T]}|\bar {\bfu}(t)|_{0,2} \le M_0$,
 such that
$ \bar{\bfu}(t) ={\bfu}_{\Omega'}(t)$ in $\bfH_\sigma(\Omega')$ for any
open and relatively compact subset $\Omega'$ of $\Omega(t,0)$ and
any $t \in [0,T]$, where $M_0$ is the same constant as in (2.8) and
${\bfu}_{\Omega'}(t)$ is the function constructed in Lemma~3.2,
corresponding to~$\Omega'$. Moreover, for $\bfu_n$ defined by~(2.11),
 $$ \bfu_n(t) \to \bar{\bfu}(t)\text{ weakly~in ~}\bfH_\sigma(\Omega),~~
   |\bar{\bfu}(x,t)| \le p(x,t)
   \text{ for a.e.~} x \in \Omega,~\forall t \in [0,T],
\eqno{(3.9)}$$
and $\bar{\bfu}=\bfu$ a.e.~on $Q$.}\vspace{0.5cm}

\noindent
{\bf Proof.} For each $t\in [0,T]$, by (1.1), the set 
$\Omega(t,0)$ is a countable union of non-decreasing, open and relatively 
compact subsets $\Omega'_i$ in $\Omega$ , $i \in {\bf N}$, such that 
$p(x,t)>\kappa_i$ on~$\Omega'_i$, $\kappa_i\to 0$.
By Lemma~3.2, the limit
$ \lim_{i\to \infty}\int_{\Omega'_i}\bfu_{\Omega'_i}(t)\cdot 
{\boldsymbol{\mathcal \xi}} dx$
exists for all ${\boldsymbol{\mathcal \xi}} \in 
\bfH_\sigma(\Omega)$.
This limit is linear and bounded with respect to 
$\boldsymbol{\mathcal \xi}$
in $\bfH_\sigma(\Omega)$ and determines a unique element $\bar{\bfu}(t)$
in ${\bfH}_\sigma(\Omega)$ by the formula
  $$ \int_\Omega \bar{\bfu}(t)\cdot \boldsymbol{\mathcal \xi} dx
   := \lim_{i\to \infty}\int_{\Omega'_i}
   \bfu_{\Omega'_i}(t)\cdot \boldsymbol{\mathcal \xi} dx,
     ~~\forall \boldsymbol{\mathcal \xi} \in \bfH_\sigma(\Omega),~~\forall t \in [0,T]. 
$$
We also have, by Lemma~3.2, $\sup_{t \in [0,T]}|\bar{\bfu}(t)|_{0,2} \le M_0$ 
and  $\bar{\bfu}(x,t)= 0$ for a.e.~$x \in \Omega$ such that
$p(x,t)=0$. Indeed, taking $\bar{\bfu}(t)$ as ${\boldsymbol{\mathcal \xi}}$
above, we get
 $$\int_\Omega |\bar{\bfu}(x,t)|^2 dx 
     =\int_{\Omega(t,0)}|\bar{\bfu}(x,t)|^2dx,$$
which implies that $|\bar{\bfu}(x,t)|=0$ for a.e.~$x \in \Omega-\Omega(t,0)
=\{x \in \Omega~|~p(x,t)=0\}$. Thus (3.9) is obtained.

Finally, we
show that $\bar{\bfu}=\bfu$ a.e.~on $Q$. To do so, for any 
$\boldsymbol{\mathcal \zeta} \in L^2(0,T;\bfH_\sigma(\Omega))$, we observe from (3.9) and Lemma 3.2 that
  $$ \int_0^T(\bar{\bfu}(t),\boldsymbol{\mathcal \zeta}(t))_\sigma dt 
      =\lim_{i \to \infty}\int_0^T\int_{\Omega'_i}
    {\bfu}_{\Omega'_i}\cdot \boldsymbol{\mathcal \zeta} dx dt
   =\int_0^T (\bfu(t),\boldsymbol{\mathcal \zeta}(t))_\sigma dt, 
  $$
which implies that $\bar{\bfu}=\bfu$ a.e.~on $Q$. Thus
$\bar{\bfu}$ is the required function. \qed
\vspace{0.5cm}

By virtue of Corollary 3.1, we may identify the function $\bfu$
with $\bar{\bfu}$; namely $\bfu$ is a function defined for every 
$t \in [0,T]$ with values in $\bfH_\sigma(\Omega)$.\vspace{0.5cm}

\noindent
{\bf Corollary 3.2.} {\it Let $\Omega'$ be an open set in $\Omega$ and 
$0\le T_1 <T'_1\le T$. Assume that $\tilde \Omega' \times [T_1,T'_1] \subset
\hat Q(p=\infty)$. Then $\bfu \in W^{1,2}(T_1,T'_1; \bfW^*_\sigma(\Omega'))$
and hence $\bfu$ is absolutely continuous as a function from
$[T_1,T'_1]$ into ${\bfW}^*_\sigma(\Omega')$.}\vspace{0.3cm}

\noindent
{\bf Proof.} Let $\bfz$ be any function in 
$C^1_0(T_1,T'_1; {\bfW}_\sigma(\Omega'))$ and take a (large) positive number
$\kappa$ and a positive integer $n(\kappa)$ so that
 $$ \|\bfz \|_{C([T_1,T'_1];{\bfW}_\sigma(\Omega'))} < \frac {\kappa}{2L_0}$$
and
 $$ p_n(x,t) > \frac {\kappa}2,
  ~~\forall (x,t) \in \Omega'\times [T_1,T'_1],~\forall n\geq n(\kappa). $$
Then we observe that
 $$ |\bfz(x,t)| \le L_0|\bfz|_{1,4} <\frac \kappa 2 < p_n(x,t),
    ~~\forall (x,t)\in \Omega'\times [T_1,T'_1],~\forall n\geq n(\kappa).$$
Therefore, just as in the proof of Lemma 3.1, we have (3.8) and
 $$ 
  \left|\int_{T_1}^{T_1'}({\boldsymbol u}_n(t), {\boldsymbol z}'(t))_\sigma dt 
   \right| =
 \left|\int_{T_1}^{T_1'}({\boldsymbol u}'_n(t), {\boldsymbol z}(t))_\sigma dt 
   \right| 
\le M_0+M_3
  \|{\boldsymbol z}\|_{L^2(T_1,T'_1;{\boldsymbol W}_\sigma(\Omega'))}.
 $$
The right hand side of the last inequality can be changed to $(M_0+M_3)
  \|{\boldsymbol z}\|_{L^2(T_1,T'_1;{\boldsymbol W}_\sigma(\Omega'))}$
with the same trick as in the proof of Lemma 3.1: whenever $\|z\|_{L^2(T_1,T'_1;{\bfW}_\sigma(\Omega'))}<1$, we~put
$$
\tilde{\bfz} = \frac{\bfz}{\|\bfz\|_{L^2(T_1,T'_1;{\bfW}_\sigma(\Omega'))}}
$$
above. Finally,
letting $n \to \infty$,
 we obtain
\begin{eqnarray*}
 \left|\int_{T_1}^{T_1'}({\boldsymbol u}(t), {\boldsymbol z}'(t))_\sigma dt 
   \right| 
&\le& (M_0+M_3)
      \|{\boldsymbol z}\|_{L^2(T_1,T'_1;{\boldsymbol W}_\sigma(\Omega'))}
 \end{eqnarray*}
for all $\bfz \in C^1_0(T_1,T'_1; \bfW_\sigma(\Omega')$. This shows that
$\bfu\in W^{1,2}(T_1,T'_1; \bfW_\sigma^*(\Omega'))$, whence 
$\bfu$ is absolutely
continuous as a function from $[T_1,T'_1]$ into $\bfW_\sigma^*(\Omega')$.
\qed
\vspace{0.5cm}

\noindent
{\bf Lemma 3.3. }{\it Assume (1.1). 
Take $\kappa>0$ and let $\Omega'$ be an open set in $\Omega$ 
such that $\tilde \Omega'\times [T_1,T_1']
 \subset \hat Q(p>\kappa)$,
where $\tilde \Omega'$ is the relative closure of 
$\Omega'$ in $\Omega$. Then, for $\bfu_n$ defined by~(2.11), we have
${\boldsymbol u}_n \to {\boldsymbol u}$ (strongly) in 
$L^2(T_1,T'_1; {\boldsymbol H}_\sigma(\Omega'))$ as $n \to \infty$.
}\vspace{0.3cm}

\noindent
{\bf Proof. } 
On account of the Aubin's compactness lemma (see [20; Lemma 5.1]), for any 
$\eps>0$ there exists a positive 
constant $A_\eps$ such that 
$$
|\bfz |_{0,2}^2\le \eps |\bfz|_{1,2}^2+ A_\eps |\bfz|_{-1,\frac 43}^2, \forall 
\bfz\in\bfV_\sigma(\Omega').
$$ 
Thus, 
$$
\int_{T_1}^{T_1'} |\bfu_n - \bfu|_{0,2}^2\,dt
\le 
\eps \int_{T_1}^{T_1'} |\bfu_n-\bfu|_{1,2}^2\,dt
+ A_\eps\int_{T_1}^{T_1'} |\bfu_n- \bfu |_{-1,\frac 43}^2\,dt.
$$
By Lemma~3.2, the last term tends to~$0$. 
Therefore, by (2.8) of Proposition 2.1 we derive from the above inequality
$$
\limsup_{n\to\infty}\int_{T_1}^{T_1'} |\bfu_n - \bfu|_{0,2}^2\,dt
\le 
\frac{2\eps}{\nu}M_0.
$$
Since $\eps$ is arbitrary, this gives the statement of the lemma.\qed
\vspace{0.5cm}

\section{The proof of Theorem 1.1}\label{sec:mainproof}

In all this section, $\bfu_n$ is the sequence defined by~(2.11).
We identify the function $\bfu$ with $\bar{\bfu}$ constructed in
Corollary 3.1; hence we have:
 $$ \bfu_n(t) \to \bfu(t)~\text{ weakly~in~}\bfH_\sigma(\Omega),~~\forall
    t \in [0,T], \eqno{(4.1)} $$
 $$ |\bfu(x,t)| \le p(x,t),~~\text{ a.e.~}x \in \Omega,~~\forall t \in
    [0,T].
     \eqno{(4.2)} $$
Furthermore, we have the following lemma.\vspace{0.3cm}

\noindent
{\bf Lemma 4.1. }
 {\it ${\boldsymbol u}_n \to {\boldsymbol u}$ in 
$L^2(0,T;{\boldsymbol H}_\sigma(\Omega))$ as $n\to\infty$.}
\vspace{0.3cm}

\noindent
{\bf Proof. }  
Let $\eps$ be any positive number and consider ${\hat Q(p\le\kappa)}$
with $\kappa:=(\frac {\eps}{18T|\Omega|})^{\frac 12}$, $|\Omega|$ being
the volume of $\Omega$.  
By (2.1), for a large integer $n_1(\eps)$ we have
$$
|\bfu_n(x,t)|\le p_n(x,t)\le p(x,t)+\kappa\le 2\kappa,~~ 
\text{ for a.e.~}(x,t)\in\hat Q(p\le\kappa),~\forall n>n_1(\eps).
$$
 Therefore, using (4.2) noted above,
$$
\int\int_{\hat Q(p\le\kappa)}|\bfu_n(x,t)-\bfu(x,t)|^2dxdt
\le
9T\,|\Omega|\,\kappa^2 = \frac{\eps}{2}.
\eqno{(4.3)} $$
Next, consider ${\hat Q(p>\kappa)}$. Take any $\kappa'\in(0,\kappa)$. 
Note that $\hat Q (p>\kappa)\subset\hat Q (p>\kappa')$ and that
by (1.1) we can find a finite number of cylindrical domains of the form 
${\Omega}_i\times [\tau_i, \tau_i']$, $i=1,2,\ldots, N$, such that 
$\tilde {\Omega}_i$ (the relative closure  of $\Omega_i$ in $\Omega$)
is contained in $\Omega$ and
$$
\hat Q (p>\kappa)
\subset  
\bigcup_{i=1}^{N}{\Omega}_i\times [\tau_i, \tau_i']
\subset 
\bigcup_{i=1}^{N}\tilde{\Omega}_i\times [\tau_i, \tau_i'] 
\subset
\hat Q (p>\kappa').
$$
Indeed, for any $ (x,t)\in\hat Q(p>\kappa)$ there exists an open set
$\Omega(x,t)\subset{\Omega}$ with
 $\tilde \Omega(x,t)\subset {\Omega}$, and there exist $\tau:=\tau(x,t), \, 
\tau':=\tau'(x,t)$ with $\tau<\tau'$ such that $t \in [\tau,\tau']$ and
$ \Omega(x,t)\times[\tau,\tau']\subset  \hat Q(p>\kappa')$. We take a finite 
covering of $\hat Q(p>\kappa)$
from this family. For such a finite covering ${\Omega}_i\times [\tau_i, \tau_i']$, $i=1,2,\ldots, N$,
it follows from Lemma~3.3 that there is a positive integer $n_2(\eps)$ such that for all $n>n_2(\eps)$
 $$
  \int_{\hat Q(p>\kappa)}|{\bfu}_n(x,t)-{\bfu}(x,t)|^2dxdt
\le
\sum_{i=1}^N \|{\bfu}_n-{\bfu}\|^2_{L^2(\tau_i,\tau_i'; {\bfH}_\sigma(\Omega_i))}\le \frac{\eps}{2}.
 \eqno{(4.4)} 
$$
Summing (4.3) and (4.4), we obtain $\int_Q |{\bfu}_n-{\bfu}|^2dxdt \le {\eps}$.
The lemma is  proved.
\qed 
\vspace{0.5cm}

We are now in a position to prove Theorem 1.1.\vspace{0.5cm}

\noindent
{\bf Proof of Theorem 1.1:}
We consider the approximate problems $P(p_n;\bfg, \bfu_{0n})$ defined by 
Definition 2.1 with $p_n$ satisfying (2.1) and 
$\bfu_{0n} = (1-\hat\delta_n)^+\bfu_0$, with $\hat\delta_n$ as in 
Proposition 2.1.
Let ${\boldsymbol v}$ be any test function from $\boldsymbol{\cal K}(p)$.
Then, for some positive constant $\delta$ we have 
${\rm supp}({\boldsymbol v}) \subset \hat Q(p > \delta)$,
so that
by Lemma 2.3 its approximate sequence
  $$ {\boldsymbol z}_n(t)=(1-\delta_n)^+ {\boldsymbol v}(t) \in K(p_n;t),
     ~~\forall t \in [0,T],~~\forall n, 
  $$
satisfies
  $$ {\rm supp}({\boldsymbol z}_n) \subset \hat Q(p>\delta),
    ~~{\boldsymbol z}_n \to {\boldsymbol v}~~\text{ in~}C^1([0,T];
     {\boldsymbol W}_\sigma(\Omega)).$$

We can now take ${\boldsymbol z}_n $ as test function in (2.4) to get
for all large $n$ with $\delta_n <1$ that
$$
\begin{array}{l}
  \displaystyle{
 ({\boldsymbol u}'_n(\tau), {\boldsymbol u}_n(\tau)
  -{\boldsymbol v}(\tau))_\sigma
  +\nu \langle {\boldsymbol u}_n(\tau), {\boldsymbol u}_n(\tau)
  -{\boldsymbol v}(\tau)
  \rangle_\sigma }\\[0.2cm]
    \displaystyle{~~~~~~~~~~~~~~~~~~~+ (\bfG (\tau,{\boldsymbol u}_n(\tau)), 
   {\boldsymbol u}_n(\tau)-{\boldsymbol v}(\tau))_\sigma} \\[0.3cm]
\displaystyle{ \le  ({\boldsymbol g}(\tau), {\boldsymbol u}_n(\tau)
-{\boldsymbol v}(\tau))_\sigma
  - \delta_nY_n(\tau)-\delta_n Z_n(\tau), ~~\text{ a.e.~}\tau \in [0,T],}
\end{array} \eqno{(4.5)} $$
where
  $$
   Y_n(\tau):= ({\boldsymbol u}'_n(\tau),{\boldsymbol v}(\tau))_\sigma $$
and
  $$ Z_n(\tau):=
    \nu \langle {\boldsymbol u}_n(\tau), {\boldsymbol v}(\tau)\rangle_\sigma 
    + (\bfG (\tau,{\boldsymbol u}_n(\tau)), {\boldsymbol v}(\tau))_\sigma 
    -({\boldsymbol g}(\tau), {\boldsymbol v}(\tau))_\sigma.
   $$
Here we note that
  $$ \left|\int_0^t Y_n(\tau) d\tau \right| \le 
   |({\boldsymbol u}_0,{\boldsymbol v}(0))_\sigma|
    + |({\boldsymbol u}_n(t),{\boldsymbol v}(t))_\sigma| 
     +T\|{\boldsymbol u}_n\|_{L^\infty(0,T;{\boldsymbol H}_\sigma(\Omega))}
      \|{\boldsymbol v}'\|_{C([0,T];{\boldsymbol H}_\sigma(\Omega))} $$
and this is uniformly bounded on $[0,T]$. Moreover, 
 \begin{eqnarray*}
 \int_0^t |Z_n(\tau)|d\tau &\le& \nu \|{\boldsymbol u}_n\|
  _{L^2(0,T;{\boldsymbol V}_\sigma)}\|{\boldsymbol v}\|
  _{L^2(0,T;{\boldsymbol V}_\sigma)}\\
 & &  + 9L_3\,\|{\boldsymbol u}_n\|_{L^2(0,T;{\boldsymbol H}_\sigma)}
       \|{\boldsymbol u}_n\|_{L^2(0,T;{\boldsymbol V}_\sigma)}
       \|{\boldsymbol v}\|_{C([0,T];{\boldsymbol W}_\sigma(\Omega))}\\[0.2cm]
 & &  +\|{\boldsymbol g}|_{L^2(0,T;{\boldsymbol H}_\sigma)}
       \|{\boldsymbol v}\|_{L^2(0,T;{\boldsymbol H}_\sigma(\Omega))},~~
  \forall t \in [0,T]
 \end{eqnarray*}
and the right hand side is bounded in 
$n$ on account of estimate (2.8).
Therefore, we have
  $$ \delta_n\int_0^t\{Y_n(\tau)+Z_n(\tau)\}d\tau \to 0
  ~~\text{ uniformly~on~}[0,T].\eqno{(4.6)} $$

After integration of (4.5) in time over $[0,t]$, 
use the integration by parts 
in the resultant and recall that 
$(G(\tau,{\boldsymbol u}_n(\tau)), {\boldsymbol u}_n(\tau)) = 0$.
Then,
$$ \begin{array}{l}
\displaystyle{\int_0^t({\boldsymbol v}', {\boldsymbol u}_n-{\boldsymbol v})
  _\sigma d\tau
   +\frac 12 |{\boldsymbol u}_n(t)-{\boldsymbol v}(t)|^2_{0,2}}\\[0.3cm]
  \displaystyle{~~~~~~~~~~~~~
     +\nu \int_0^t\langle {\boldsymbol u}_n, {\boldsymbol u}_n
    -{\boldsymbol v}\rangle_\sigma d\tau 
       - \int_0^t \langle \bfG(\tau,{\boldsymbol u}_n),
     {\boldsymbol v}\rangle_\sigma d\tau }\\[0.3cm]
  \displaystyle{ \le  \int_0^t({\boldsymbol g}, {\boldsymbol u}_n
    -{\boldsymbol v})_\sigma  d\tau   + \frac 12 |{\boldsymbol u}_n(0)-{\boldsymbol v}(0)|^2_{0,2}
  - \delta_n\int_0^t\{Y_n+Z_n\}d\tau.}
\end{array} \eqno{(4.7)} $$
Now we pass to the limit $n\to\infty$ in (4.7).
The first term of the left hand side and the others are bounded from below
by the respective terms with the limit $\bfu$, since ${\boldsymbol u}_n \to
{\boldsymbol u}$ weakly in $L^2(0,T;{\boldsymbol V}_\sigma(\Omega))$ and 
${\boldsymbol u}_n(t) \to{\boldsymbol u}(t)$ weakly in ${\boldsymbol H}_\sigma
(\Omega)$ for every $ t \in [0,T]$ by (4.1). This also allows to pass 
to the limit in the first term of the right hand side. The initial condition is chosen so that 
$\bfu_n(0)=\bfu_{0n} \to \bfu_0$ in~$\bfW_\sigma(\Omega)$. As for the nonlinear term, we observe that
\begin{multline*}
\int_0^t (\bfG(\tau,{\boldsymbol u}_n) - (\bfG(\tau,{\boldsymbol u}(\tau)),
{\boldsymbol v})_\sigma d\tau\\
= \sum_{k,j=1}^3 \int_0^t\int_\Omega 
\left\lbrace 
v^{(j)}\left(u^{(k)}_n - u^{(k)}\right)\frac{\partial u_n^{(j)}}{\partial x_k}
+
v^{(j)}u^{(k)}\frac{\partial \left(u^{(j)}_n - u^{(j)}\right)}{\partial x_k}
\right\rbrace dx\,d\tau
\to 0
\end{multline*}
by virtue of strong convergence ${\boldsymbol u}_n \to {\boldsymbol u}$ in 
$L^2(0,T; \bfH_\sigma(\Omega))$ (cf.  Lemma~4.1),
combined with the weak convergence in $L^2(0,T; \bfV_\sigma(\Omega))$. 
Consequently, for any~${\boldsymbol v} \in \boldsymbol{\cal K}(p)$ it follows from 
(4.6) and (4.7) that
\begin{multline*}
\int_0^t({\boldsymbol v}', {\boldsymbol u}-{\boldsymbol v})_\sigma d\tau
 +\frac 12 |{\boldsymbol u}(t)-{\boldsymbol v}(t)|^2_\sigma
    +\nu \int_0^t \langle {\boldsymbol u}, {\boldsymbol u}- {\boldsymbol v} 
\rangle_\sigma d\tau \\
 + \int_0^t \int_\Omega ({\boldsymbol u}\cdot\nabla){\boldsymbol u}\cdot
      \nabla({\boldsymbol u}-{\boldsymbol v})dx d\tau 
     \\
\le \int_0^t ({\boldsymbol g},{\boldsymbol u}-{\boldsymbol v})_\sigma d\tau
     + \frac 12 |{\boldsymbol u}_0-{\boldsymbol v}(0)|^2_\sigma,
\end{multline*}
i.e.~$\bfu$ solves~(1.2) . With (4.2), the condition (ii) of Definition~1.1 is satisfied.

It remains to show (i). 
As for the initial condition, $\bfu(0) = \bfu_0$ 
since  $\bfu_{0n}$ converges strongly to $\bfu_0$ 
in $\bfW_\sigma(\Omega)$. 
Let us show that $t\mapsto (\bfu(t), \bfv(t))_\sigma$ is of bounded 
variation on $[0,T]$ for every $\bfv \in {\boldsymbol{\cal K}}(p)$.
For simplicity, we denote
   $$  m_{\bfv}=\|\bfv(t)\|_{L^\infty(0,T;{\bfH}_\sigma(\Omega))},~~~
    M_{\bfv}=\|\bfv(t)\|_{C([0,T];{\bfW}_\sigma(\Omega))}.$$
Since
supp$(\bfv) \subset \hat Q(p>\kappa)$ for a certain $\kappa>0$, 
we can find a finite covering
$ \Omega_i\times [T_i,T'_i],~~i=1,2,\cdots, N$, of supp$(\bfv)$ such that
  $$ {\rm supp}({\bfv}) \subset \bigcup_{i=1}^N \Omega_i\times [T_i,T'_i]
     \subset \bigcup_{i=1}^N \tilde \Omega_i\times [T_i,T'_i]
     \subset \hat Q(p>\kappa).$$
For each $i$, it follows from Lemma 3.1 that the restriction of $\bfu$ to
$[T_i,T'_i] \times \Omega_i$ is of bounded variation as a function from
$[T_i,T'_i]$ into $\bfW^*_\sigma(\Omega_i)$; we denote by $V_i(\bfu)$ its total
variation.
Now, for any $s, t \in [T_i,T'_i]$ we observe that
  \begin{eqnarray*} 
 |(\bfu(t),\bfv(t))_\sigma -(\bfu(s),\bfv(s))_\sigma|
&\le& |(\bfu(t)-\bfu(s),\bfv(t))_\sigma + (\bfu(s),\bfv(t)-\bfv(s))_\sigma|\\
 &\le& M_{\bfv}|\bfu(t)-\bfu(s)|_{{\bfW}_\sigma^*(\Omega_i)}
      +m_{\bfu}|\bfv(t)-\bfv(s)|_{{\bfH}_\sigma(\Omega)}.
\end{eqnarray*}
The total variation of $t \mapsto (\bfu(t), \bfv(t))_\sigma$ on the interval
$[T_i,T'_i]$ is bounded by 
$M_{\bfv} V_i(\bfu)+m_{\bfu}\int_{T_i}^{T'_i}|\bfv'|_{0,2}dt$. Therefore, the
total variation on the whole interval $[0,T]$ is not larger than
  $$ M_{\bfv}\sum_{i=1}^N V_i(\bfu)
       +m_{\bfu}\int_0^T|\bfv'|_{0,2}dt ~(<\infty).$$
The proof of Theorem 1.1 is now complete. \hfill $\Box$
\vspace{0.5cm}

\noindent
{\bf Remark 4.1.} In particular, if the support of the test function 
$\bfv$ is contained in $\hat Q(p=\infty)$, 
then it follows from Corollary 3.1 that
the function $t \mapsto (\bfu(t),\bfv(t))_\sigma$ is absolutely continuous on
$[0,T]$ and $\bfu(x,0)=\bfu_0(x)$ on $\{x \in \Omega~|~p(x,0)=\infty\}$ pointwise.
\vspace{0.3cm}

\noindent
{\bf Remark 4.2.} Assume that $\bfg\equiv 0$, 
i.e.~no external forces are present. In this case, if
at time $t_0$ the whole region $\Omega$ 
is blocked by a total obstacle, 
i.e.~$p(x,t_0)=0$ for all  $x\in \Omega$, then
the flow vanishes starting from the moment $t_0$: 
$\bfv\equiv 0$ in $[t_0,T]\times \Omega$. 
In other words, if the obstacle grows to the whole region
at some time $t_0$,
then it blocks the flow efficiently, 
even if the obstacle itself diminishes afterwards. 
In fact, since 
$\bfv(\cdot,t_0)\equiv 0$ in $\Omega$, it follows that $\bfv\equiv 0$ is
the trivial solution of the Navier--Stokes equation on 
$(t_0,T)\times \Omega$ and this is the solution of the variational 
inequality of the Navier--Stokes type, which can be constructed in our
approximate procedure, too.
\vspace{1cm}

\begin{center}
{\bf References}
\end{center}
\begin{enumerate}
\item H. Abels, Longtime behavior of solutions of a Navier--Stokes/Cahn--Hillard
system, pp. 9--19 in {\it Nonlocal and abstract parabolic equations and their
applications}, Banach Center Publ. {\bf Vol. 86}, Polish Acad. Sci., Inst. 
Math., Warsaw, 2009.
\item W. H. Alt and I. Pawlow, Existence of solutions for non-isothermal phase
separation, Adv. Math. Sci. Appl., {\bf 1} (1992), 319--409.
\item M. Biroli, Sur l’in\'equation d'{\'e}volution de Navier--Stokes.
C. R. Acad. Sci. Paris Ser. A-B, 275 (1972), A365--A367.
\item M. Biroli, Sur in\'equation d'{\'e}volution de Navier--Stokes.
Nota I, Atti Accad. Naz. Lincei Rend. Cl. Sci. Fis. Mat. Natur. (8), 52
(1972), 457--460.
\item M. Biroli, Sur in\'equation d'{\'e}volution de Navier--Stokes.
Nota II. Atti Accad. Naz. Lincei Rend. Cl. Sci. Fis. Mat. Natur. (8), 52
(1972), 591--598.
\item M. Biroli, Sur in\'equation d'{\'e}volution de Navier--Stokes.
Nota III, Atti Accad. Naz. Lincei Rend. Cl. Sci. Fis. Mat. Natur. (8), 52 (1972), 811--820
\item H. Br\'ezis, {\it Op\'eratuers Maximaux Monotones et 
Semi-groupes de 
Contractions dans les espaces de Hilbert}, Math.~Studies 5, North-Holland,
Amsterdam, 1973.
\item M. Brokate and J. Sprekels, {\it Hysteresis and Phase Transitions}, {\bf Vol. 121}, Springer-Verlag, Berlin--Heidelberg--New York, 
1996.
\item P. Colli, N. Kenmochi and M. Kubo, A phase-field model with temperature
dependent constraint, J. Math. Anal. Appl., {\bf 256} (2001), 668--685.
\item A. Damlamian, Some results on the multi-phase Stefan problem, Comm.
Partial Differential Equations, {\bf 2} (1977), 1017--1044.
\item  H.J. Eberl, D.F. Parker and M.C.M. van Loosdrecht, 
A new deterministic
spatio-temporal continuum model for biofilm development, 
Computational and Mathematical Methods in Medicine, {\bf 3} (2001), 161--175.
\item L. C. Evans and R. F. Gariepy, {\it Measure Theory and Fine Properties of Functions }, CRC Press, Boca Raton--London--New York--Washington, D.C., 1992.
\item T. Fukao and N. Kenmochi, Variational inequality for the Navier-Stokes
equations with time-dependent constraint, Gakuto Internat. Ser. Math. Sci. Appl., {Vol. 34} (2011), 87--102.
\item T. Fukao and N. Kenmochi, Quasi-variational inequalities 
approach to 
heat convection problems with temperature dependent velocity constraint, §Discrete Contin. Dyn. Syst., 
{\bf 35} (2015), 2523--2538.
\item D. Gilbarg and N. S. Trudinger, {\it Elliptic Partial Differential 
Equations of Second Order}, Springer-Verlag, Berlin--Heidelberg--New York, 1998.
\item H. Inoue and M. \^Otani, Periodic problems for heat convection 
equations in
noncylindrical domains, Funkcail. Ekvac., {\bf 40} (1997), 19--39.
\item N. Kenmochi, Solvability of nonlinear evolution equations with 
time-dependent constraints and applications, Bull. Fac. Edu., Chiba Univ., {\bf 30} (1981), 1--87.
\item N. Kenmochi and M. Niezg\'odka, Viscosity approach to modelling
non-isothermal diffusive phase separation, Jpn. J. Ind. Appl. Math.,
{\bf 13} (1996), 135--169.
\item M. Kubo, A. Ito and N. Kenmochi, Non-isothermal phase separation models:
weak well-poseness and global estimates, Gakuto Internat. Ser. Math. Sci. Appl.,
{\bf Vol. 14} (2000), 311--323. 
\item J. L. Lions, {\it Quelques m\'ethodes de r\'esolution des probl\`emes aux 
limites non lin\'eaires}, Dunod Gauthier--Villars, Paris, 1969.
\item Y. Murase and A. Ito, Mathematical model for the process of brewing 
Japanese sake and its analysis, Adv. Math. Sci. Appl. {\bf 23} (2013), 297--317.
\item M. \^Otani, Nonmonotone perturbations for nonlinear parabolic equations
associated with subdifferential operators, J. Differential Equations, 
{\bf 46} (1982), 
268--299.
\item  M. Peszy\'nska, A. Trykozko, G, Iltis and S. Schlueter, 
Biofilm growth in
porous media: Experiments, computational modeling at the porescale, 
and upscaling, Advances in Water Resources, 1--14, 2015.
\item G. Prouse, On an inequality related to the motion, in
any dimension, of viscous, incompressible fluids. Note I, Atti Accad. Naz. Lincei Rend. Cl. Sci. Fis. Mat. Natur. (8), 67 (1979), 191--196.
\item G. Prouse, On an inequality related to the motion,
in any dimension, of viscous, incompressible fluids. Note II, Atti Accad. Naz. Lincei Rend. Cl. Sci. Fis. Mat. Natur. (8), 67 (1979), 282--288.
\item K. Shirakawa, A. Ito, N. Yamazaki and N. Kenmochi, Asymptotic stability
for evolution equations governed by subdifferentials, pp. 287-310 in {\it
Recent Development in Domain Decomposition Methods and Flow Problems}, Gakuto Internat. Ser. Math. Sci. Appl., {\bf Vol. 11}, Gakk\=otosho, Tokyo, 1998.
\item R. Temam, {\it Navier-Stokes Equations, Theory and Numerical Analysis},
North-Holland, Amsterdam, 1984.
\item Y. Yamada, On nonlinear evolution equations generated by the 
subdifferentials, J.~Fac. Sci. Univ. Tokyo, Sect. IA, {\bf 23}(1976), 491--515.
\item K. Yosida, {\it Functional Analysis} (Sixth edition), Springer-Verlag,
Berlin--Heidelberg--New York, 1980.

\end{enumerate}

\end{document}